\documentclass[11pt]{article}
\usepackage{fourier}
\usepackage[scaled=0.85]{berasans}
\usepackage[scaled=0.85]{beramono}
\usepackage{amsmath}
\usepackage{url}
\newcommand\bLP{\\[\bigskipamount]}
\newcommand\ZZ{\mathbb{Z}}
\newcommand\de\delta
\newcommand\om\omega
\newcommand\iy\infty
\newcommand\thalf{\tfrac12}
\newcommand{\hyp}[5]{\,\mbox{}_{#1}F_{#2}\!\left(
  \genfrac{}{}{0pt}{}{#3}{#4};#5\right)}

\begin{document}
\title{On an identity by Chaundy and Bullard. II. More history}
\author{Tom H. Koornwinder and Michael J. Schlosser\footnote{Partly
supported by FWF Austrian Science Fund grant S9607.}}
\date{}
\maketitle
\begin{abstract}
An identity by Chaundy and Bullard writes $1/(1-x)^n$ ($n=1,2,\ldots$)
as a sum of two truncated binomial series. In a paper which appeared in 2008
in Indag.\ Math.\ the authors surveyed many aspects of this identity.
In the present paper we discuss much earlier occurrences of this identity
in works by Hering (1868), de Moivre (1738) and de Montmort (1713).
A relationship with Krawtchouk polynomials in work by Greville (1966) is also
discussed.
\end{abstract}
%
\section{Introduction}
In our paper \cite{1} we surveyed the history of the often rediscovered formula
\begin{equation}
1=(1-x)^{n+1}\sum_{k=0}^m\binom{n+k}k x^k
+x^{m+1}\sum_{k=0}^n\binom{m+k}k (1-x)^k.
\label{1}
\end{equation}
We attributed the formula to Chaundy \& Bullard \cite[p.256]{2} (1960).
However, we later learnt that some giant steps back in time can be made to
much earlier occurrences of this formula.
Almost one century before Chaundy \& Bullard
the formula was given by Hering \cite{3} (1868). Then, with a jump of more than
one century, the formula was found in the work of de Moivre \cite{6} (1738).
Even 25 years earlier the formula was given in implicit form already by
de Montmort \cite{8} (1713).

The paper successively discusses these three early occurences of the formulas.
Next a correspondence between Samuel Pepys and Isaac Newton, having some relation
with identity \eqref{1}, is briefly discussed.
We conclude with a much more recent connection with Krawtchouk polynomials
which is implicit in Greville \cite{13} (1966).
\bLP
{\bf Acknowledgements}\quad
We thank Pieter de Jong for communicating us the occurrence of formula
\eqref{1} in de Moivre \cite{6}, \cite{7}, for calling our attention
to reference \cite{18}, and for sending us his manuscript \cite{5}.
We also thank the Mathematics Department of Lund University for
sending us a
copy of Hering's paper \cite {3}, which happened to be in the Small Boxes of
the G\"osta Mittag-Leffler Separate Collection, stored at this
Department (see \cite{18}).
\section{Hering (1868)}
In 1868 Hering \cite[p.14, formula 1)]{3} derived:
\begin{equation}
(1-x)_n^{-m}=(1-x)^{-m}-(1-x)^{-m}\,x^n\,(1-\overline{1-x}\,)_m^{-n}.
\label{7}
\end{equation}
Here $(1-x)_n^{-m}$ is the power series in $x$ of $(1-x)^{-m}$ cut after the
$n$-th term. Similarly, $(1-\overline{1-x}\,)_m^{-n}$ is the power series
in $1-x$ of $(1-(1-x))^{-n}$ cut after the $m$-th term. Thus
Hering already had \eqref{1}.

Hering's proof is different from any of the proofs given in \eqref{1}.
For generic non-integer $m$ he writes for the left-hand side of \eqref{7}:
\begin{equation}
(1-x)_n^{-m}=\sum_{k=0}^{n-1}\frac{(m)_k}{k!}\,x^k
=\frac{(m)_{n-1}}{(n-1)!}\,x^{n-1}\,\hyp21{-n+1,1}{-m-n+2}{x^{-1}}
=\frac{(m)_{n-1}}{(n-1)!}\,\frac{x^n}{x-1}\,\hyp21{-m+1,1}{-m-n+2}{\frac1{1-x}}\,,
\label{8}
\end{equation}
where inversion of the order of summation is used in the second equality and
Pfaff's transformation formula in the third equality.
If $m$ tends to a positive integer, the last
${}_2F_1$ becomes
\begin{equation}
\sum_{k=0}^{m-1}\frac{(-m+1)_k}{(-m-n+2)_k}\,(1-x)^{-k}
+\sum_{k=m+n-1}^\iy\frac{(-m-n+k+2)_{n-1}}{(-m-n+2)_{n-1}}\,(1-x)^{-k}.
\label{9}
\end{equation}
(Here, although not emphasized by Hering, we should require for convergence
that $|x-1|>1$. This can later be relaxed in \eqref{7} by analytic continuation.)$\;$
After multiplication of \eqref{9} by $\binom{m+n-2}{n-1}x^n/(x-1)$ we can rewrite the first term
(by inversion of the order of summation) as
\[
-(1-x)^{-m}\,x^n\,\sum_{k=0}^{m-1}\frac{(n)_k}{k!}\,(1-x)^k=
-(1-x)^{-m}\,x^n\,(1-\overline{1-x}\,)_m^{-n},
\]
and the second term as
\[
(-1)^n(1-x)^{-m-n}x^n\,\sum_{k=0}^\iy\frac{(n)_k}{k!}\,(1-x)^{-k}
=(1-x)^{-m}.
\]
Thus by substitution in \eqref{8} Hering settled \eqref{7}.

Formula \eqref{7} is just one of many formulas derived in \cite{3}.
Hering does not specially emphasize this particular result.
\section{de Moivre (1738)}
A much earlier reference was kindly communicated to us by Pieter de Jong
and also mentioned in his manuscript \cite{5}.
In 1738 A.~de Moivre \cite[p.196]{6} (see also the 1754 edition
\cite[p.224]{7}) wrote:
\begin{quotation}
{\sl But as there is a particular elegancy for the Sums of a finite number of
Terms in those Series whose Coefficients are figurate numbers beginning at Unity,
I shall Set down the Canon for those Sums.

Let $n$ denote the number of Terms whose Sum is to be found, and $p$ the
rank or order which those figurate numbers obtain, then the Sum will be}
\[
\frac{1-x^n}{(1-x)^p}-\frac{nx^n}{(1-x)^{p-1}}
-\frac{n\,.\,n+1\,.\,x^n}{1\,.\,2\,.\,(1-x)^{p-2}}
-\frac{n\,.\,n+1\,.\,n+2\,.\,x^n}{1\,.\,2\,.\,3\,.\,(1-x)^{p-3}}
-\frac{n\,.\,n+1\,.\,n+2\,.\,n+3\,.}{1\,.\,2\,.\,3\,.\,4\,.\,(1-x)^{p-4}}\,,\;
{\rm etc.}
\]
{\sl which is to be continued till the number of Terms  be $=p$.}
\\{}
[In the numerator of the last term above a factor $x^n$ is missing. This may have
been a printer's error.]
\end{quotation}
According to de Moivre \cite[Corollary at end of p.195]{6}
the figurate numbers of order $p$ are the successive coefficients
in the power series of $\tfrac1{(1-x)^p}$. Thus, since he begins at unity,
these are binomial coefficients $\binom{p+k-1}k$ ($k=0,1,2,\ldots$).
This is slightly different from the modern definition
given by Dickson \cite[p.7]{4}, who defines
the $k$-th figurate number of order $p$
as the binomial coefficient $\binom{p+k-1}p$.
The difference is that de Moivre starts counting orders at 1
(for instance triangular numbers have order 3 for him),
while Dickson starts counting them at 0, by which triangular numbers have
order 2.

Thus by the above quotation de Moivre
gives the identity
\[
\sum_{k=0}^{n-1} \binom{p+k-1}{k} x^k
=(1-x)^{-p}-(1-x)^{-p}\,x^n\sum_{k=0}^{p-1}\binom{n+k-1}k (1-x)^k.
\]
Indeed, this shows that de Moivre had already \eqref{1} in 1738.

The section ``Of the Summation of recurring Series'' starting
in de Moivre \cite[p.193]{6} gives some indication how he
obtained his result. We will summarize this in modern terminology and we
state everything at once for general $p$ instead of stating it for
$p=1,2,3$, etc.

First de Moivre discusses infinite power series
$S=\sum_{k=0}^\iy c_k x^k$ in which the coefficients satisfy a recurrence
relation
$c_k=a_1c_{k-1}+a_2c_{k-2}+\cdots+a_p c_{k-p}$ with coefficients
$a_j$ independent of $k$. Then he observes that
$S=q(x)/(1-a_1x-\cdots-a_px^p)$, where $q(x)$ is a polynomial of degree at
most $p-1$ which can be explicitly computed.
Next he observes that the figurate numbers
\[
c_k=\binom{p+k-1}k=\frac{(k+1)_{p-1}}{(p-1)!}\,,
\]
being polynomials of degree $p-1$ in $k$ which vanish for
$k=-1,-2,\ldots,-p+1$,
are annihilated by the $p$-th finite difference:
\[
\sum_{l=0}^p(-1)^l\,\binom pl\,c_{k-l}=0\qquad(k=1,2,\ldots).
\]
From this he derives that $S=(1-x)^{-p}$ in this case.
Finally he applies the same method to terminating power series
$S_n=\sum_{k=0}^{n-1} c_k x^k$. A recurrence relation
$c_k=a_1c_{k-1}+a_2c_{k-2}+\cdots+a_p c_{k-p}$ will then yield
$S_n=(q(x)+x^n\,r(x))/(1-a_1x-\cdots-a_px^p)$ for certain polynomials
$q(x)$ and $r(x)$ of degree at most $p-1$. For $c_k$ being
the figurate numbers
we get $S_n=(1+x^nr(x))/(1-x)^p$, from which $r(x)$ can be computed,
in principle.
However, de Moivre does not give an argument how he
arrives at the nice explicit expression for $r(x)$ expanded in powers of $1-x$.
Probably, he found the expression for low values of $p$ and then extrapolated.

De Moivre also made an important step by which he might have concluded
the multi-variable generalization
of \eqref{1} given in \cite[(1.5)]{1} and first obtained (as far as we know)
by Damjanovic, Klamkin and Ruehr \cite{17} in 1986.
In fact, de Moivre \cite[Problem LXIX, p.191]{6},
\cite[pp.50,51]{7} gave for any of the $n$ summands in the outer sum in
\cite[(1.5)]{1} a probabilistic interpretation coming from the problem of
points with $n$ players (see the case of 2 players below).
Adding up these chances to 1 would have given him the multi-variable formula,
just as was pointed out in \cite{17} (in detail in a situation with three urns
by Bosch and Steutel). We are puzzled why de Moivre missed this final step,
and also why he did not give a probabilistic interpretation of \eqref{1}.
\section{de Montmort (1713)}
In 1713 appeared the second edition of the {\em Essay d'analyse sur les jeux de
hazard} \cite{8} by Pierre Raymond de Montmort. It contained among others
a new solution of the so-called {\em problem of points} for two players.
This problem comes from  a game of chance with two players Pierre and Paul who
have chances
$p$ and $1-p$, respectively, of winning each round. The player who
has first won a certain number of rounds (this number may be different
for Pierre and Paul) will collect the entire prize.
Suppose that the game is prematurely interrupted when Pierre has to win
still $n$ rounds and Paul $m$ rounds. What is then a fair division of the stake?
See Hald \cite[\S14.1]{10} for a description how this problem was handled by
de Montmort.

In the case of equal chances the problem was already solved by Pascal and
Fermat in 1654.
In the case of unequal chances Johann Bernoulli generalized their solution.
Bernoulli gives his solution in a letter to de Montmort
dated 17 March 1710. This letter is included in the second edition of
de Montmort's book, see \cite[pp.283--298]{8}, in particular p.295
(English translation available at \cite{9}).
De Montmort also gives this solution in his main text, see
\cite[pp. 244--245, \S190]{8} (English translation at~\cite{9}).
Curiously, Bernoulli is not mentioned there by de Montmort.
Neither he acknowledges this new result of Bernoulli in his polemical
discussion of earlier work on the problem of points
in the Avertissement of the second edition of his book. This discussion
starts on p.~xxxiv of \cite{8} (English translation available at~\cite{9}).
Bellhouse \cite{16} gives an interesting discussion of the relationship between
de Montmort and de Moivre.

Bernoulli's solution is as follows.
Imagine Pierre and Paul still play $m+n-1$ rounds. Then there will certainly be a
winner. Pierre will be the winner if he has won $n$ or more of these rounds.
His chance for this
is
\begin{equation}
\sum_{k=0}^{m-1}\binom{m+n-1}k p^{m+n-k-1}(1-p)^k.
\label{2}
\end{equation}
But if Pierre has won $n-1$ or less of these rounds then Paul will be the winner.
The chance for this is
\begin{equation}
\sum_{\ell=0}^{n-1} \binom{m+n-1}\ell
p^\ell(1-p)^{m+n-\ell-1}.
\label{3}
\end{equation}
The two chances add up to 1, both by the probabilistic
interpretation and by the binomial formula.

De Montmort \cite[p.245, \S191]{8} continues to give other expressions for the chances
for Pierre respectively Paul to win.
Imagine they still play until there is a winner.
Pierre will be the winner if he has won already $n-1$ rounds and Paul at most
$m-1$ rounds, and if then the next round is won by Pierre.
Thus the chance for Pierre to win is
\begin{equation}
p\,\sum_{k=0}^{m-1} \binom{n+k-1}k\,p^{n-1}(1-p)^k,
\label{4}
\end{equation}
Similarly, for Paul the chance to win is\begin{equation}
(1-p)\,\sum_{k=0}^{n-1} \binom{m+k-1}k\,p^k(1-p)^{m-1}.
\label{5}
\end{equation}
These two chances necessarily add up to 1.
Thus the resulting formula (not given by de Montmort) is
\[
p^n\,\sum_{k=0}^{m-1} \binom{n+k-1}k\,(1-p)^k+
(1-p)^m\,\sum_{k=0}^{n-1} \binom{m+k-1}k\,p^k=1.
\]
by which formula \eqref{1} is proved in a probabilistic way.
This is essentially the same proof as was quoted from much more recent
literature in \cite[end of Section 6]{1}.

Clearly the chances \eqref{2} and \eqref{4} are the same.
This is not explicitly observed by de Montmort, but it is indicated in the example
where $n=5$ and $m=3$.
In the general case
the resulting identity is \cite[(2.7)]{1}. There we referred to Guenther \cite{11},
who gave various proofs and references  (but none older than
1933) for this identity, including the
probabilistic proof we just observed.
\section{Pepys and Newton (1693)}
In 1693 Samuel Pepys wrote a letter to Isaac Newton with a question about
a probabilistic problem coming from a question to Pepys by John Smith
(see \cite{2}, \cite{19}). The question was (in modern terms):
\begin{quote}
Let $6k$ fair dice be tossed independently and suppose that at least $k$ "6"'s appear. For which $k=1,2,3$
this has the greatest chance to happen?
\end{quote}
Newton wrote back three times. He answered correctly that the case $k=1$ has
the highest probability. He actually computed the probability for $k=1$ and 2.
He also gave a theoretical argument about which Stigler \cite{18}, as late
as 2006, observed that it was incorrect.
Chaundy \& Bullard \cite{2} showed more generally:
\begin{quote}
We work with fair dice with $s$ faces. Let $g(sn,n)$ be the chance that a
selected face turns up less than $n$ times in $sn$ throws. Then $g(sn,n)$
increases with $n$ for fixed $s$.
\end{quote}
They proved this statement by expressing $g(sn,n)$ in terms of \eqref{2} and
then using that \eqref{2} is equal to \eqref{4}.
(In passing, in connection with their proof of this identity, they observed
the identity \eqref{1}.)\;
Finally, by working with \eqref{4} they could prove their claim quoted above.
\section{Greville (1966)}
This last item does not push the history of identity \eqref{1} further back, but mentions
an unexpected aspect of this identity which is offered by Greville
\cite[p.166]{13}. A further description is given in \cite[\S4,6]{14}.
Greville considers the smoothing filter $f\to g$ given by
\begin{equation*}
g(y)=\sum_{x=-N}^N f(y-x)\,{\bf K}_{2n}(x,0)\,w(x)\qquad(y\in\ZZ),
\end{equation*}
where $w(x):=\binom{2N}{N+x}$ and
${\bf K}_n$ is the Christoffel-Darboux kernel for the
orthogonal polynomials $p_n$ satisfying
\begin{equation*}
\sum_{x=-N}^N p_n(x)\,p_m(x)\,w(x)=h_n\,\de_{n,m}\qquad(n,m\in\{0,1,\ldots,2N\}).
\end{equation*}
Then the polynomials $p_n$ are special shifted Krawtchouk polynomials
\begin{equation*}
p_n(x)=K_n(x+N;\thalf;2N),
\end{equation*}
but this is not explicitly mentioned by Greville. Then we also see that
$h_n=2^{2N}\binom{2N}n^{-1}$ and that
\begin{equation*}
{\bf K}_{2n}(x,0)=\sum_{k=0}^{2n} \frac{p_k(x)p_k(0)}{h_k}\,.
\end{equation*}
Greville wants to compute the characteristic function (or transfer function)
$\phi$ associated with
this smoothing filter, given by
\begin{equation*}
\phi(\om):=\sum_{x=-N}^N {\bf K}_{2n}(x,0) e^{-i\om x}.
\end{equation*}
Then he derives that
\begin{equation}
\phi(\om)=1-(\sin^2(\om/2))^{n+1}P(\sin^2(\om/2))=
(\cos^2(\om/2))^{N-n}Q(\sin^2(\om/2))
\label{6}
\end{equation}
for certain polynomials $P$ of degree $N-n-1$ and $Q$ of degree $n$.
Then, with the same argument as in \cite[Section 6.1]{15} and \cite[Remark 2.2]{1},
Greville explicitly obtains $P$ and $Q$. As a consequence, \eqref{6}
takes the form of \eqref{1} with $m=N-n-1$. Greville also concludes from
the explicit expression that $\phi$ is monotonically decreasing from 1 to 0 on $[0,\pi]$.
Later Herrmann \cite{12} independently computed \eqref{6} in a different way
in order to arrive at this result of monotonical decrease of $\phi$, which he called
maximal flatness.

\quad\\
\begin{footnotesize}
\begin{quote}
T. H. Koornwinder, Korteweg-de Vries Institute, University of
Amsterdam,\\
P.O.\ Box 94248, 1090 GE Amsterdam, The Netherlands;\\
email: {\tt T.H.Koornwinder@uva.nl}
\bLP
M.~J. Schlosser,
Fakult\"at f\"ur Mathematik, Universit\"at Wien,\\
Nordbergstrasse 15,
A-1090 Vienna, Austria;\\
email: {\tt michael.schlosser@univie.ac.at}
\end{quote}
\end{footnotesize}


\begin{thebibliography}{99}
%
\bibitem{16}
D. Bellhouse.
{\em Banishing Fortuna: Montmort and De Moivre},
J. History Ideas  69 (2008),  559--581.
%
\bibitem{2}
T. W. Chaundy and J. E. Bullard,
{\em John Smith's problem}, Math. Gazette 44 (1960), 253--260.
%
\bibitem{18}
T. Claesson and J. Peetre,
{\em Index of the G\"osta Mittag-Leffler separate collection. Part I, Small Boxes},
Lund University, 1996;
\url{http://staff.science.uva.nl/~thk/art/2012/ChaundyBullard2/}
%
\bibitem{15}
I. Daubechies,
{\em Ten Lectures on Wavelets}, Regional Conference Series in Applied Math.,
vol. 61, SIAM, 1992.
%
\bibitem{4}
L. E. Dickson,
{\em History of the theory of numbers. Vol. II: Diophantine analysis},
Carnegie Inst. of Washington, Washington, D.C., 1920; reprinted, Chelsea,
New York, 1966.
%
\bibitem{14}
E. Diekema and T. H. Koornwinder, {\em Differentiation by integration using
orthogonal polynomials, a survey}, J. Approx. Theory 164 (2012), 637--667;
{\tt arXiv:1102.5219v2 [math.CA]}. 
%
\bibitem{13}
T. N. E. Greville,
{\em On stability of linear smoothing formulas},
SIAM J. Numer. Anal. 3 (1966), 157--170.
%
\bibitem{11}
W. C. Guenther,
{\em Solution of Problem E 1829},
Amer. Math. Monthly 74 (1967), 1134--1135.
%
\bibitem{10}
A. Hald,
{\em A history of probability and statistics and their applications before 1750},
Wiley, New York, 1990.
%
\bibitem{3}
A. G. Hering,
{\em Summation der $n$ ersten Glieder der binomischen Reihe mittelst der 
Theorie der hypergeometrischen Reihen},
Programm der Realschule in Chemnitz, 1868;
\url{http://staff.science.uva.nl/~thk/art/2012/ChaundyBullard2/};
JFM 01.0089.04.
%
\bibitem{12}
O. Herrmann,
{\em On the approximation problem in nonrecursive digital filter design},
IEEE Trans. Circuit Theory 18 (1971), 411--413.
%
\bibitem{5}
P. de Jong, {\em The arrangement of the arithmetical triangle}, manuscript, 2012.
%
\bibitem{1}
T. H. Koornwinder and M. J. Schlosser,
{\em On an identity by Chaundy and Bullard. I},
Indag. Math. (N.S.) 19 (2008), 239--261; {\tt arXiv:0712.2125v3 [math.CA]}. 
%
\bibitem{6}
A. de Moivre,
{\em The doctrine of chances}, Second edition, Woodfall, London, 1738;
reprinted,
Frank Cass \& Co., Ltd., London, 1967.
%
\bibitem{7}
A. de Moivre,
{\em The doctrine of chances}, Third edition, A. Millar, London,
1756; reprinted, Chelsea, New York, 1967.
%
\bibitem{8}
P. R. de Montmort,
{\em Essay d'analyse sur les jeux de hazard},
Seconde \'edition, Laurent le Conte, Paris, 1713; reprinted,
Chelsea, New York, 1980.
%
\bibitem{9}
R. J. Pulskamp,
{\em Pierre Raymond de Montmort},
\url{http://www.cs.xu.edu/math/Sources/Montmort/montmort.html}.
%
\bibitem{19}
S. M. Stigler,
{\em Isaac Newton as a probabilist},
Statist. Sci. 21 (2006), 400--403.
%
\bibitem{17}
{\em Solution of Problem 85-10}, SIAM Rev. 28 (1986) 243--244.
%
\end{thebibliography}
\end{document}